\documentclass [12pt,a4paper]{article}
\usepackage{amssymb, theorem}
\def\proof{\noindent{\it Proof.} }

\def\bbbr{{\mathbb R}}
\def\bbbc{{\mathbb C}}
\def\bbbz{{\mathbb Z}}
\def\Ker{\mathrm{Ker}\,}
\def\Rng{\mathrm{Rng}\,}

\def\Tr{\mathrm{Tr}\,}

\def\diag{\mathrm{Diag}}

\def\iH{{\cal H}}

\def\iF{{\cal F}}

\def\ot{\otimes}
\def\osum{\oplus}
\def\det{\mathrm{Det}\,}

\def\framedformula#1#2{$$\vcenter{\hrule\hbox{\vrule\kern5pt
\vbox{\kern5pt\hbox{$\displaystyle#1$}%
\kern5pt}\kern5pt\vrule}\hrule}\eqno#2$$}

\newtheorem{thm}{Theorem}[section]
\newtheorem{prop}{Proposition}[section]
\theorembodyfont{\rm}

\newtheorem{pl}{Example}[section]

\def\qed{\nobreak\hfill $\square$}
\begin{document}
\rightline{Acta Math. Hungar., to appear}
\ \vskip 1cm 
\centerline{\LARGE {\bf Strongly subadditive functions}}
\medskip
\bigskip
\bigskip
\centerline{\bf Koenraad Audenaert$^1$, Fumio Hiai$^2$ and D\'enes Petz$^3$}
\begin{center}
$^1$ Department of Mathematics, Royal Holloway, University of London,\\
Egham, Surrey TW20 0EX, UK

\medskip
$^2$ Graduate School of Information Sciences,
Tohoku University, \\ Aoba-ku, Sendai 980-8579, Japan 

\medskip
$^3$ Alfr\'ed R\'enyi Institute of Mathematics, \\
H-1364 Budapest, POB 127, Hungary
\end{center}

\bigskip
\begin{abstract}
Let $f:\bbbr^+ \to \bbbr$. The subject is the trace inequality 
$\Tr f(A) + \Tr f(P_2AP_2) \le \Tr f(P_{12}AP_{12}) + \Tr f(P_{23}AP_{23})$, 
where $A$ is a positive operator, $P_1,P_2,P_3$ are orthogonal projections 
such that $P_1+P_2+P_3=I$, $P_{12}=P_1+P_2$ and  $P_{23}=P_2+P_3$. There 
are several examples of functions $f$ satisfying the inequality 
(called (SSA)) and the case of equality is described. 

 

\end{abstract}

\section{Introduction}

Matrix monotone and matrix concave functions play important roles in several
applications. Assume that $f:\bbbr^+\to \bbbr$ is a continuous function. It
is {\it matrix monotone} if $0 \le A \le B$ implies $f(A)\le f(B)$ for every
matrix $A$ and $B$. The function $f$ is called {\it matrix concave} if one of the 
following two equivalent conditions holds:
\begin{equation}
f(\lambda A +(1-\lambda) B) \ge \lambda f(A) + (1-\lambda)f(B)
\end{equation}
for every number $0 < \lambda < 1$ and for positive definite square matrices
$A$ and $B$ (of the same size). In the other condition the number $\lambda$ 
is (heuristically) replaced by a matrix:
\begin{equation}
f(CAC^* +DBD^*) \ge C f(A) C^* + Df(B)D^*
\end{equation}
if $CC^*+DD^*=I$, see the books \cite{Bh, PD} about the details. It is
surprising that a matrix monotone function is matrix concave.

Motivated by some applications we study the functions $f$ which are
strongly subadditivite in the following sense. Let $P_1,P_2,P_3$ be 
orthogonal projections such that $P_1+P_2+P_3=I$. Then
\begin{equation}\label{E:ssa1}
\Tr f(A) + \Tr f(P_2AP_2) \le \Tr f(P_{12}AP_{12}) + \Tr f(P_{23}AP_{23}),
\end{equation}
where $P_{12}:=P_1+P_2$ and  $P_{23}:=P_2+P_3$. The special case when $P_2=0$
could be called {\it subadditivity}. This holds for any concave function
\cite[Theorem 2.4]{HP2}.

The first example $f(x)=\log x$ appeared already \cite{AP}, here we have several 
other examples and a sufficient condition.  The strongly subadditive functions
are concave in the sense of real variable and all known examples are matrix
concave.

\section{Motivation}
The {\it second quantization} in quantum theory is mathematically a procedure
which associates an operator on the Fock space $\iF(\iH)$ to an operator
on the Hilbert space $\iH$ \cite{B-R}. The simplest example is $\iH=\bbbc$, 
then $\iF(\iH)$ is $\ell^2(\bbbz^+)$. To the number $\mu>0$ (considered 
as a positive operator) we associate $\Gamma (\mu)$ defined as
$$
\Gamma (\mu) \delta_n= \mu^n \delta_n \qquad (n=0,1,2,\dots),
$$
where $\delta_n$ are the standard basis vectors. $\Gamma$ can be extended
to arbitrary finite dimension by the formula
$$
\Gamma (H_1 \osum H_2)=\Gamma (H_1) \ot \Gamma (H_2).
$$
In this way to any positive operator $H \in B(\iH)$ we have a positive
operator $\Gamma (H) \in B(\iF(\iH))$. The construction of a statistical
operator, analogue of the Gaussian distribution, is slightly more complicated.
For a positive operator $A$ set
$$
\alpha(A)=\frac{\Gamma(H)}{\Tr \Gamma(H)},\quad \mbox{where}\quad H=A (I+A)^{-1}.
$$
In particular, if $A=\lambda$, then
$$
\alpha (\lambda) \delta_n=
\frac{1}{1+\lambda}\left(\frac{\lambda}{1+\lambda}\right)^n \delta_n\,.
$$
The von Neumann entropy 
$$
S(\alpha (A)):=- \Tr \alpha (A) \log \alpha (A)
$$
equals to $\Tr \kappa(A)$, where $\kappa(x):=-x \log x + (x+1)\log (x+1)$
\cite{B-R, JPP}.

From the formula
$$
\log x = \int_0^\infty \frac{1}{1+t}-\frac{1}{x+t}\, dt\,.
$$ 
we get
\begin{eqnarray*}
\kappa(x)&=& -\int_0^\infty \frac{x}{1+t}-\frac{x}{x+t}\,dt
+\int_0^\infty \frac{x+1}{1+t}-\frac{x+1}{x+1+t}\,dt\cr &=&
\int_0^1 1-\frac{x}{1+t}-\frac{t}{x+t}\,dt +\int_1^\infty \frac{x+1}{t}-\frac{x}{1+t}
-\frac{1}{x+t}\,dt.
\end{eqnarray*}
Since both integrands are matrix concave, the integrals are matrix concave, too.

$$
\kappa'(x)=\log \Big(1+\frac{1}{x}\Big) >0
$$
and $\kappa$ is monotone. Hence $\kappa(x) \ge \kappa(0) = 0$. The positivity 
together with matrix concavity implies matrix monotonicity, \cite{HP}.

Let $\iH=\iH_1\oplus \iH_2 \oplus \iH_3$ be a finite dimensional Hilbert space
and let
$$
A=\left[\matrix{ A_{11} & A_{12}& A_{13} \cr A_{12}^* & A_{22} & A_{23} \cr 
A_{13}^*  & A_{23}^* & A_{33} }\right],
$$
be a positive invertible operator and  set
$$
B=\left[\matrix{ A_{11} & A_{12} \cr A_{12}^* & A_{22} }\right], \qquad
C=\left[\matrix{ A_{22} & A_{23} \cr A_{23}^* & A_{33} }\right].
$$

The strong subadditivity of the von Neumann entropy,
$$
S( \alpha(A)) + S( \alpha(A_{22})) \le S( \alpha(B)) + S( \alpha(C)),
$$
has the equivalent form
\begin{equation}\label{E:kappa}
\Tr \kappa(A) + \Tr \kappa(A_{22}) \le \Tr \kappa(B) + \Tr \kappa(C).
\end{equation}
The case of equality is studied in the paper \cite{JPP} and the general
properties of entropy are in the book \cite{PD}.

\begin{prop} \label{P:JPP}
The equality
$$
\Tr \kappa(A) + \Tr \kappa(A_{22}) \le \Tr \kappa(B) + \Tr \kappa(C),
$$
in the strong subadditivity holds if and only if $A$ has the form
\begin{equation}\label{E:trivi}
A=\left[\matrix{ A_{11} &
\left[\matrix{a & 0}\right]&
0
\cr 
&&\cr
\left[\matrix{a^* \cr 0}\right]&
\left[\matrix{c & 0 \cr 0 & d}\right]&
\left[\matrix{0 \cr b }\right] \cr &&\cr
 0 &
\left[\matrix{0 & b^*}\right]&
A_{33}
 }\right]=\left[\matrix{
\left[\matrix{ A_{11} & a \cr a^* & c}\right] &
0 \cr & \cr & \cr
0 &
\left[\matrix{d & b \cr b^* & A_{33}}\right]
 }\right], 
\end{equation}
where the parameters $a, b, c, d$ (and $0$) are matrices.
\end{prop} 

Note that the matrix $c$ or $d$ in the theorem can be $0 \times 0$.

We are interested in the (differentiable) functions $f$ such that the inequality
\framedformula{\Tr f(A) + \Tr f(A_{22}) \le \Tr f(B) + \Tr f(C)}{(SSA)}
holds. We call this {\it strong subadditivity} for the function $f$.
The strong subadditivity holds for the function $\kappa$. Another equivalent
formulation of the strong subadditivity is (\ref{E:ssa1}).

\section{Particular examples}

\begin{pl}
If
\begin{equation}\label{abcd}
A=\left[\matrix{ a & 0 & d \cr 0 & b & 0 \cr 
d^*  & 0 & c }\right]
\end{equation}
is a numerical matrix, then it is an exercise to show that (SSA) holds for
this kind of $A$ if and only if $f$ is a concave function.
\end{pl}

\begin{pl}
The strong subadditivity does not hold for the function $f(t)=-1/t$.
The following counterexample is due to Ando \cite{Ando}:
Let
$$
X \equiv A^{-1}:=
\left[\matrix{ 4 & 8 & -2 \cr 8 & 20 & 0 \cr -2 & 0 & 9 }\right].
$$
Then 
$$
A=\left[\matrix{ \frac{45}{16} & -\frac{9}{8} & \frac{5}{8} 
\cr  -\frac{9}{8} &  \frac{1}{2} &  -\frac{1}{4} 
\cr \frac{5}{8} & -\frac{1}{4} & \frac{1}{4} }\right].
$$
We have
$$
\Tr A^{-1}=33, \quad \Tr A_{22}^{-1}=2, \quad \Tr B^{-1}= \frac{212}{9}, 
\quad \Tr C^{-1}=12 
$$
and (SSA) becomes
$$
33 + 2 \ge \frac{212}{9}+ 12
$$
and this is not true.
\end{pl}

\begin{pl}
It is elementary that the strong subadditivity holds for the function $f(t)=-t^2$.
The equality holds if and only $A_{13}=0$. \qed
\end{pl}

\begin{pl}\label{P:log}
It was proved in \cite{AP} that the strong subadditivity holds for the function 
$f(t)=\log t$ and the equality holds if and only if $A_{13}=A_{12}A_{22}^{-1}A_{23}$.

We present an alternative approach. Now (SSA) is equivalent to
$$
\det A\cdot\det A_{22}\le\det B\cdot\det C.
$$
Let
$$
\hat A:=\diag(A_{11}^{-1/2},A_{22}^{-1/2},A_{33}^{-1/2})
A\,\diag(A_{11}^{-1/2},A_{22}^{-1/2},A_{33}^{-1/2}).
$$
Then (SSA) is equivalent to 
$$
\det\hat A\le\det\hat B\cdot\det\hat C.
$$ 
In other words, we may assume that the diagonal of $A$ consists of $I$'s.
Since
$$
\left[\matrix{ I&-\hat A_{12}&0\cr 0&I&0\cr 0&0&I}\right]
\left[\matrix{I&\hat A_{12}&\hat A_{13}\cr \hat A_{12}^*&I&\hat A_{23}^*\cr
\hat A_{13}^*&\hat A_{23}^*&I}\right]
\left[\matrix{ I&0&0\cr 0&I&-\hat A_{23}\cr 0&0&I}\right]\qquad \qquad
$$
$$\qquad \qquad
=\left[\matrix{I-\hat A_{12}\hat A_{12}^*&0&\hat A_{13}-\hat A_{12}\hat A_{23}
\cr \hat A_{12}^*&I&0\cr \hat A_{13}^*&\hat A_{23}^*&I-\hat A_{23}^*\hat A_{23}
}\right],
$$
equality holds in (SSA) if $\hat A_{13}=\hat A_{12}\hat A_{23}$, equivalently
$A_{13}=A_{12}A_{22}^{-1}A_{23}$.  This condition is sufficient for the equality. \qed
\end{pl}

\begin{pl}
Since
$$
{1\over2}\Big(A+\diag(1,1,-1)A\,\diag(1,1,-1)\Big)
=\left[\matrix{ A_{11}&A_{12}&0\cr A_{12}^*&A_{22}&0\cr 0&0&A_{33}}\right]
$$
we get a majorization
$$
A\succ \left[\matrix{ A_{11}&A_{12}&0\cr A_{12}^*&A_{22}&0\cr 0&0&A_{33}}
\right],
$$
that is, the eigenvalue vector $\vec\lambda(A)$ majorizes that of
$$
\left[\matrix{ A_{11}&A_{12}&0\cr A_{12}^*&A_{22}&0\cr 0&0&A_{33}}\right].
$$
For any concave  function $f$, this implies that $f\circ\vec\lambda(A)$ is 
weakly majorized by the $f\circ\lambda$ of
$$
\left[\matrix{ A_{11}&A_{12}&0\cr A_{12}^*&A_{22}&0\cr 0&0&A_{33}}\right]
$$ 
so that
\begin{equation}\label{acd}
\Tr f(A)\le\Tr f\left(\left[\matrix{ A_{11}&A_{12}&0\cr
A_{12}^*&A_{22}&0\cr 0&0&A_{33}}\right]\right)
=\Tr f(B)+\Tr f(A_{33}).
\end{equation}
Hence
$$
\Tr f(A)+\Tr f(A_{22})\le\Tr f(B)+\Tr f(A_{22})+\Tr f(A_{33}).
$$
This says that if $A_{23}=0$ (or $A_{12}=0$), then (SSA) holds for every 
concave function $f$.

Note that inequality (\ref{acd}) is written as 
$$
\Tr f(A) =  \Tr Pf(A)P+ \Tr Qf(A)Q \le \Tr f(PAP + QAQ) 
$$
when $P$ and $Q$ are orthogonal projections and $P+Q=I$. This is a special case of
Jensen's trace inequality for concave functions \cite[Theorem 2.4]{HP2}. \qed
\end{pl}

\begin{pl}

The representation
\begin{equation}\label{E:power}
y^t=\frac{\sin \pi t}{\pi}\int_0^\infty \frac{\lambda^{t-1}y}{\lambda+y}\,d\lambda
\end{equation}
is used to show that $f(x)=x^t$ is operator monotone when $0 < t <1$. From
this we obtain
$$
\int_0^x y^{t-1}\,dy =
\frac{\sin \pi t}{\pi}\int_0^x \int_0^\infty \frac{\lambda^{t-1}}{\lambda+y}
\,d\lambda\,dy
$$ 
which gives
$$
x^t=\frac{t \sin \pi t}{\pi}\int_0^\infty \lambda^{t-1}
(\log (x+\lambda)-\log \lambda)\,d\lambda.
$$
So we have a similar formula to (\ref{E:power}):
\begin{equation}\label{E:power2}
x^t=\frac{t \sin \pi t}{\pi}\int_0^\infty \lambda^{t-1} \log\Big(1+\frac{x}{\lambda}
\Big)\,d\lambda.
\end{equation}

Since inequality (SSA) is true for the functions $f_\lambda(x):=
\log\Big(1+\frac{x}{\lambda}\Big)$, by integration it follows for 
$x^t$ when $0 < t <1$. 

We analyze the condition for equality and use the decomposition
$\iH=\iH_1\oplus \iH_2 \oplus \iH_3$. For $f_\lambda$ the equality condition is
$$
A_{13}=A_{12}(\lambda+A_{22})^{-1}A_{23},
$$
see Example \ref{P:log}. This holds for every $\lambda>0$. If $\lambda \to \infty$ in
the relation
$$
\lambda A_{13}=A_{12}\Big[\lambda(\lambda+A_{22})^{-1}\Big]A_{23},
$$
then we conclude $ A_{13}=0=A_{12}A_{23}$. The latter condition means that $\Rng A_{23} 
\subset \Ker A_{12}$, or equivalently $(\Ker A_{12})^\perp \subset \Ker A_{23}^*$.

The linear combinations of the functions $x \mapsto 1/(\lambda+x)$ form an
algebra and due to the Stone-Weiersrass theorem $A_{12}g(A_{22})A_{23}=0$ for any 
continuous function $g$.

We want to show that the equality implies the structure (\ref{E:trivi}) of the
operator $A$. We have $A_{23} : \iH_3 \to \iH_2$ and $A_{12}:\iH_2 \to \iH_1$. 
To show  the structure (\ref{E:trivi}), we have to find a subspace $H \subset \iH_2$
such that
$$
A_{22} H \subset H, \quad H^\perp \subset \Ker A_{12}, \quad H \subset \Ker A_{32},
$$
or alternatively $( H^\perp=)K \subset \iH_2$ should be an invariant subspace of $A_{22}$
such that 
$$
\Rng  A_{23} \subset K \subset \Ker A_{12}.
$$

Let
$$
K:=\Big\{\sum_i A_{22}^{n_i}A_{23}x_i \,:\, x_i \in  \iH_3, n_i \in \bbbz^+ \Big\}
$$
be a set of finite sums.  It is a subspace of $\iH_2$. The property  $\Rng  A_{23} \subset K$
and the invariance under $A_{22}$ are obvious.  Since
$$
A_{12} A_{22}^{n}A_{23}x=0,
$$
$ K \subset \Ker A_{12}$ also follows. \qed
\end{pl}

\section{Sufficient condition}

\begin{thm}
Let $f:(0,+\infty) \to \bbbr$ be a function such that $-f'$ is matrix monotone.
Then the inequality (SSA) holds.
\end{thm}

\proof
The idea of the previous example is followed. A matrix monotone function has
the representation
$$
a+ bx+\int_0^\infty \left(\frac{\lambda}{\lambda^2+1}-
\frac{1}{\lambda+x}\right)\,d\mu(\lambda), 
$$
where $b \ge 0$, see (V.49) in \cite{Bh}. Therefore, we have the representation
$$
f(t)=c-\int_1^t \left(
a+ bx+\int_0^\infty \left(\frac{\lambda}{\lambda^2+1}-
\frac{1}{\lambda+x}\right)\,d\mu(\lambda)
\right)\,dx.
$$
By integration we have
$$
f(t)=d-at-\frac{b}{2}t^2+ \int_0^\infty 
 \left(\frac{\lambda}{\lambda^2+1}(1-t)+\log
\Big(\frac{\lambda}{\lambda+1}+\frac{t}{\lambda+1}\Big)\right)\,d\mu(\lambda).
$$
The first quadratic part satisfies the (SSA) and we have to check the integral.
Since $\log x$ is a strongly subadditive function, so is the integrand.
The integration keeps the property. \qed

The previous theorem covers all known examples, but we can get new examples.

\begin{pl}
By differentiation we can see that $f(x)=-(x+t) \log (x+t)$ with $t\ge 0$ satisfies
(SSA). Similarly, $f(x)=-x^t$ satisfies (SSA) if $1 \le t \le 2$.

In some applications \cite{PSz} the operator monotone functions
$$
f_p(x)={p(1-p)}\,{ (x-1)^2 \over (x^p -1)
(x^{1-p}-1)}\qquad (0 < p <1)
$$
appear. 

For $p=1/2$ this is an (SSA) function. Up to a constant factor, the function is
$$
(\sqrt{x}+1)^2=x+2\sqrt{x}+1
$$
and all terms are known to be (SSA).  The function $-f_{1/2}'$ is evidently 
matrix monotone.

Numerical computation shows that $-f_p'$ seems to be matrix monotone. \qed
\end{pl}

{\bf Acknowledgements.} This work was partially supported by the Hungarian 
Research Grant OTKA T068258 (D.P.) and Grant-in-Aid for Scientific Research 
(B)17340043 (F.H.) as well as by Hungary-Japan HAS-JSPS Joint Project (D.P.\ 
\& F.H.). D.P.\  is also grateful to Professor Tsuyoshi Ando for communication.
and for his example.

\end{document}